\documentclass{amsart}

\usepackage{amsmath, amsthm, amssymb, graphicx, enumerate, everypage, datetime, getfiledate, lmodern}

\newcounter{forwardtheorem}

\newcounter{citedtheorems}

\newtheorem*{theorem-x}{Theorem}
\newtheorem*{theorem-m}{Main Theorem}
\newtheorem*{theorem-n}{Main Theorem}
\newtheorem*{cor-x}{Corollary}
\newtheorem*{lemma-x}{Lemma}
\newtheorem*{concl-x}{Conclusion}
\newtheorem*{claim-x}{Claim}

\newtheorem{thm-lit}[citedtheorems]{Theorem}
\newtheorem{defn-a}[citedtheorems]{Definition}

\newcommand{\mcp}{\mathcal{P}}

\newcommand{\vp}{\varphi}

\newcommand{\br}{\vspace{4mm}}

\newcommand{\tcb}{\textcolor{blue}}

\title{Realizing infinity}

\author{Maryanthe Malliaris and Assaf Peretz}

\date{Version of February 2, 2021; minor edits July 5, 2022.}

\begin{document}

\begin{abstract}
What happens when mathematics realizes infinity. When are mathematical definitions actually useful? 
\end{abstract}

\maketitle

\br
\section*{I}
\br

Many people find reading fiction beneficial to their life. 
This is depicted in “The Lion, The Witch And The Wardrobe” where we are shown how visiting a possibly fictitious land behind the wardrobe tremendously improves the visitors’ lives. One can ask why are these things useful in people’s lives. One can also hypothesize whether the very distinction between reality and fiction is closer to a fictitious one than to a real one; more than likely a question of scale rather than any absolute distinction. 

\br
\section*{II}
\br

Mathematics has a unique place in our lives. The place of complete precision -- at least as far as we can see. This precision is powerful. In mathematics we don’t describe some triangle in the sand, but an idealized triangle, with lines which have no width and can go on infinitely long, and where the meeting place of two lines is a single point. We completely define what “a line” is 
and we (at least mathematicians) believe that these idealized ideas are useful to us in our everyday lives. We know there is no 
\textbf{real} line in the world. Yet introducing the concept, the definition, of a line has been extremely useful in our history.

Similarly introducing the concept/definition of a variable x has been extremely powerful in our history. 
{The variable x does not exist} {per se}, 
it is not real, and yet, its method has been powerful. The introduction of the ‘variable’ is a methodology, but this methodology necessitates  
a unique creation -- that of the mysterious ‘variable’.  

Already we can see that there is a difference between the introduction of a line and that of a variable. Both may be ``fictitious'' -- in the sense of not existing in the ‘real world’ -- but they are so in quite a different manner. 

Another important “fictitious,” or should we say imaginary element was defined as ``$i = \sqrt{-1}$.'' 
This too proved to be quite useful, though it took thousands of years to come to fruition, from its first appearance to the Greeks, to being useful to Euler and Gauss.

Mathematical definitions give us the power of creation. The moment we define something it exists. We breathe life into it. Contrary to most works of fiction though, mathematics, demanding complete precision, has provided itself with the quality of dis-creation, the power to annihilate from existence. If Montaigne can say that he may contradict himself but never contradict the truth, mathematics can’t do that. It does not allow contradiction, and if a definition leads to contradiction it is dis-created, and annihilated from existence. As a metaphor, we could say that math constantly creates matter, and then anti-matter, and then uses the energy created by their meeting in order to get to where it wants to get to – the goal, the Proof. 

\br
\section*{III}
\br

Infinity as a methodology existed in different manners in mathematics from early on. A major creation by Leibniz and Newton was that of the methodology of the limit which goes to infinity, a clear mathematical definition of a progression to infinity. The invention of calculus, we can say, has been useful.\footnote{Another example of a version of realization of infinity is the definition of the point of compactification. This too has proved quite useful.}

At the end of the 19th century mathematics took another step and created infinity itself. It gave a definition for, and thus created and brought to life, “infinity.” Today, more than a century later, this creation is still ignored and barely understood by all but a few mathematicians (and by most scientists, who have a presumption of 
the finitude of the physical world).   

\br
\section*{IIII} 
\br

{A leap of faith:} 

For religious people, we are finitary beings and to reach God, to reach the infinite, we need to actually cross infinity. It is an impossibility, and thus what it needs is a leap of faith. 
Until the mathematical invention of the infinite, infinity was a goal: a line can go on and on and on, numbers can strive to get there in their limit. Through this limit we could integrate and produce many useful things, but it was always an impossible goal. 
{\emph{\underline{Its realization by the use of a mere definition is what is so}}}
\\ \noindent \emph{\underline{staggering about mathematics and its power.}}

\vspace{1mm}
The moment infinity is defined, it is created and starts to have skeleton structure we can investigate. Its investigation yielded some surprising results, the first of which is perhaps Cantor realizing at the end of the 19th century that there are different sizes of infinity.

\br
\section*{IIIII}
\br

The power of creation, of realization in mathematics, is a great power, but like works of fiction, some realizations are more impactful than others. Every mathematician can invent/create a unique world with the snap of their fingers, a world which they may extensively work to describe, through proofs.
There is an area however where the creation\tcb{-}power 
of mathematics is far greater than that of fiction, and that is in its power of annihilation; of dis-creation. Unlike fiction (which has the power to harm but not to dis-create\footnote{Even God, per the bible, only destroys, not dis-creates.}) years of work can be annihilated in a minute by another mathematician’s disproof of a created world.

\br
\section*{IIIIII}
\br

Just so there isn’t the risk of being misunderstood we add the following: 

We are not saying that mathematics is fiction, nor that certain mathematical definitions are more real than others. In fact, 
Kant saw mathematics as (synthetic) a priori, that is, more basic in our existence than anything in the world. Descartes saw math as the hidden language for the book of the world. Our point is about the wondrous out-worldly power of creation by definition in mathematics.

In the history of philosophy from Plato to Kant, in different ways, mathematics had a meta-physical status, that is, the status of being beyond the physical world, however we might understand this beyond. Mathematics, in a way, is distinguished from the real, and yet we can use it and learn from it in the real world. That, in the end, is what we are doing in life in the real world.

\br
\section*{IIIIIII}
\br

{We give a mundane example of what we discussed above. 
\begin{quotation}\noindent \emph{John walked near his favorite cafe when he noticed a poolah on the side of the road.  
 We define a Poolah as an animal with three legs, no hands, three eyes, a single beak, and two tails. John hadn’t seen a poolah in a while, so he quickly proceeded towards it to take a closer look, when all of a sudden the poolah jumped and started running away.} 
\end{quotation}
Fiction can create a poolah. That is what it does. But in-real-life we can’t just create a poolah by defining it. Life doesn’t work that way.  
But mathematics does. Mathematics possesses the power of creation by definition. The moment of a mathematical definition is a moment of creation. Defining $i=\sqrt{-1}$ gives $i$ a life of its own. For example, we can derive that $i^2=-1$. 
The moment of definition lets the creation live its life -- 
or be disproved as impossible and dis-created. }

\br
\section*{IIIIIIII}
\br

$0$, $x$, $i$, infinity, limit structures, transcendentals, $\pi$, triangles, and any other mathematical definition – are some more real than others? 
Is $1+1=2$ more real than $10^{20}+10^{20}=2\times 10^{20}$ or $x+x = 2\times x$? While some people want to cling to a supposed real and fictitious mathematics, we see it as a matter of what is actually useful (which, as mentioned above, may take centuries to know).

Before delving more in-depth into these questions, and the power of the definition of the infinite, through concrete examples, we propose 
the following formal definition of the usefulness of definitions: 

\begin{defn-a} \label{d:defn} \emph{ }
\begin{enumerate}
\item  We say a definition is \emph{somewhat useful} if it is an integral part of a correct proof of a theorem.
\item We say a definition is \emph{actually useful} if it is an integral part of a correct proof of a theorem which doesn’t include the definition either in its assumptions or in its result. 
\end{enumerate}
\end{defn-a}

\section*{IIIIIIIII}
\br

Paul Cohen, in lecture notes explaining his celebrated work on the continuum hypothesis \cite{cohen}, writes that ``Cantor's work in set theory was the subject of much criticism to the effect that it dealt with fictions."

Why was this felt to be different from the use of infinity in classical mathematics? Is such a feeling of difference justified? 

It may seem at first glance that in the earlier cases, infinity's role is mainly as an organizing principle, helping us to see in a clearer and 
clearer way how finite objects behave.  
Consider statements and expressions like: 
 \emph{There is no largest prime. } 
 \emph{The sequence $1/2^n$ converges to 0 as $n \rightarrow \infty$.}
 $\int^5_3 3x^2 dx$.
  \emph{ There is a nondecreasing function $f: \mathbb{N} \rightarrow \mathbb{N}$ such that for every n, and 
    every graph G of size at least $f(n)$, $G$ contains a homogeneous set of size $n$.} 
Indeed, for many practical applications of calculus, once one knows the limit theory works, 
it is sufficient to get close enough. 

Still, this calm feeling of a steady progression to the very large, 
with steadily better information all the way up, can already be challenged in interesting cases. 
As an example from combinatorics, Szemer\'edi's regularity lemma says roughly that \emph{huge} finite graphs 
can be well-approximated by much smaller random graphs. 
That is, given $\epsilon > 0$, one defines a pair of vertex sets $(A,B)$ in a graph to be ``$\epsilon$-regular" if, informally, the edges between 
$A$ and $B$ are uniformly distributed (up to error $\epsilon$), similarly to a random bipartite graph with edge probability equal to the 
density $d(A,B)$.  
Then the lemma says for every $\epsilon > 0$ there is $N = N(\epsilon)$ such the vertices of any sufficiently large graph can be partitioned into no more than $N$ pieces of approximately equal size so that most pairs of pieces are $\epsilon$-regular.  In principle, this can give us 
a reasonable approximate picture of a truly huge finite graph. However, by a theorem of Gowers, 
$N$ necessarily has tower-type dependence on $\epsilon$, so the scale at which the picture appears is already well beyond our own. 

Such a lemma may grant us a kind of second sight, showing how certain mathematical objects behave at huge, 
but still finite, scales, which were otherwise not within our power to observe.  
Yet the behavior of graphs at the intermediate scales may remain quite mysterious. 
  
Already with infinity as a limit, 
when the vastness of the infinite takes up a certain residence among finite objects, 
the path from the small to the huge is not necessarily experienced as an unbroken path.

\br
\section*{IIIIIIIIII}
\br

To isolate what is surprising about the definition of infinity, we may try to separate what is \emph{generally} or 
\emph{often} surprising about mathematical definitions, from what may be specifically surprising in this case.

``\emph{Call the set H of vertices of the graph G homogeneous if either 
   no two vertices in H are connected by an edge, or every two distinct 
   vertices in H are connected by an edge.}''  
 ``\emph{Let $p$ be the largest prime.}'' ``\emph{Consider a geometry satisfying Euclid's first four postulates.}''
 ``\emph{Call a real number transcendental if it is not the root of any polynomial 
   with rational coefficients.}''
   
Recall that is meaningful to give a definition before we have evidence of anything satisfying it, 
and indeed to prove examples exist on general grounds without having any in hand.
If a definition eventually 
leads to a contradiction it is thrown out -- as we have said, dis-created (``let p be the largest 
prime'') -- but such a definition was not necessarily useless; proof by contradiction 
has a long history, and its inherent mysteries need not reflect a mystification coming from the infinite.  
Recall that a definition may appear to generalize when what it really does is simplify.  
Finally, a definition which encompasses many obvious or natural examples -- even, seemingly, ``minimally so'' -- doesn't a priori rule out 
less familiar examples of the same definition.

\br
\section*{IIIIIIIIIII}

\br

Let us consider\footnote{looking, we should emphasize, only as logicians of our own time, not historians of mathematics.} 
how deeply the definition of the infinite draws from the finite, beginning with Cantor's precise ideas which unify and 
organize our understanding of finite sets, their size, their extensions, and their arithmetic. 
We might say that a correct understanding of the finite opens a door by which the infinite also enters. 

First is an organization of the idea of size. 
According to Cantor's definition, two sets A, B have the same size, in symbols $|A| = |B|$, if they can be put in bijection (one-to-one correspondence). Say $|A| \leq |B|$  if A has the same size as a subset of B, and $|A| < |B|$ if $|A|\leq|B|$ but there is no bijection. These definitions make sense even if we do not know the sizes of A or B, and clearly are meaningful for finite sets. 
According to Dedekind's definition, a set is infinite if it can be put in bijection 
with a strict subset of itself.  

Notice this tells us we cannot increase the size of an infinite set by adding a single element. 

Cantor noticed that there is another operation which increases the size of any set, regardless of its size: the power set operation (recall 
that the power set of $X$, is the set of subsets of $X$). If 
$X = \{ 0, 1, 2 \}$, then 
$\mcp(X) = \{ \{\}$, $\{0\}$, $\{1\}$, $\{2\}$, $\{0,1\}$, $\{0,2\}$, $\{1,2\}$, $\{0,1,2\} \}$ and so $3 = |X| < |\mcp(X)| = 8$. 
According to Cantor's famous theorem, 
$|X| < |\mcp(X)|$ for any set $X$, thus, there are different infinite sizes, and what's more, no largest infinite size. 

Meanwhile, suppose we write $0 = \emptyset$, $1 = \{ 0 \}$, $2 = \{ 0, 1 \}$ and in general ${n = \{ 0, \dots, n-1 \}}$.   
Then in von Neumann's language, each natural number $n$ is a set which is \emph{transitive} and \emph{linearly ordered by $\in$}; 
call such a set an ordinal. They can be thought of as order-types of well ordered sets. 
The class of ordinals are linearly ordered, and closed under successor (if $\alpha$ is an ordinal, so is 
$\alpha \cup \{ \alpha \}$, called ``$\alpha + 1$'') and under union. 
So the natural numbers are not the only ordinals: we have their 
union $\omega = \{ 0, 1, 2, \cdots \}$, then $\omega + 1$, and so on.  
Define the cardinals to be ordinals which cannot be put in bijection with any smaller 
ordinal. 
Then any ordinal $\alpha$ has a unique size, or \emph{cardinality}, $|\alpha|$ -- the cardinal of the same size.  
The cardinals inherit a well-ordering from the ordinals. 
(Assuming the well-ordering principle, any set can be well-ordered, then is order-isomorphic to a unique ordinal, and so has a well-defined 
size or ``{cardinality}'' too.) 

The event of defining infinite sizes/cardinals makes the infinite literal and useful. 
Cantor's cardinals are linearly ordered: 
\[	0, 1, 2, 3, ... \aleph_0, \aleph_1, \aleph_2, ... \]
There are arithmetic operations which 
specialize to the usual operations when the input is finite.  However, once one of the cardinals involved is infinite, 
new phenomena appear. 
For instance, the 
fundamental theorem of cardinal arithmetic says:
\begin{quotation}
	Let $\kappa$, $\lambda$ be nonzero cardinals, at least one of them infinite.
	\\ Then $\kappa + \lambda = \kappa \times \lambda = \max \{ \kappa, \lambda \}$.
\end{quotation}

\br
\section*{IIIIIIIIIIII} 

\br

We may test Cantor's discovery against Definition \ref{d:defn}, taking as evidence one of his own proofs. 

As above, call a real number algebraic if it is the root of some polynomial with rational coefficients, and transcendental if it is not. It had been an open question to demonstrate the existence of transcendental numbers. Using his new understanding of infinite sizes, Cantor solved this problem by the following proof. Call a set \emph{countable} if it can be put in bijection with the natural numbers.  Each polynomial in x with rational coefficients is determined by a finite sequence of rational numbers, whose length depends on its degree. The rationals are countable.  For fixed n, the set of length n sequences of elements of a countable set is countable. The union of countably many at most countable sets is countable. So the set of finite sequences of rationals is countable. Each of our countably many polynomials has at most finitely many roots. And thus the number of algebraic numbers is countable. The reals have the same size as the power set of $\mathbb{N}$, and so are uncountable. Thus, not only do transcendental numbers exist, there are many more of them than algebraic numbers. Note that this proof does not produce any examples of transcendental numbers. 

This proof allegedly caused a sensation. 

Moreover, Cantor's cardinals revealed ``hidden in plain sight'' hard questions about objects people thought they knew well. 
For example, 
$|X| < |P(X)|$ begs the question: how much larger? Taking $X = \{ 0, 1, 2 \}$, 3 is quite a bit smaller than 8.  Since the reals have the same size as the power set of the natural numbers, this begs the question of whether there is an infinite set of reals ``of intermediate size,'' too large to be put in bijection with $\mathbb{N}$ and not large enough to be put in bijection with 
$\mathbb{R}$. This led, of course, to the continuum hypothesis, Hilbert's first problem (on his list of open questions presented at the 1900 International Congress of Mathematicians).  
Addition and multiplication of infinite cardinals worked so simply, it's interesting that exponentiation appears more difficult. 

Thus already from here we can see that, according to Definition \ref{d:defn}, Cantor's definition is \emph{actually useful}. 

\br

\section*{IIIIIIIIIIIII}

\br

Insofar as infinite cardinals and ordinals come ``out from under the overcoat'' of the finite, we might expect that they would carry 
certain useful features of the natural numbers into new areas and new domains.  An example of this is that proof by 
induction can now be extended to transfinite induction, and definition by recursion to transfinite recursion.

For example, a theorem in the wonderful book of Komjath and Totik \cite{kt} says that it is possible to decompose Euclidean space, $\mathbb{R}^3$, as the disjoint union of circles  
 of radius 1. 
 The proof begins by enumerating the points of $\mathbb{R}^3$ as $\langle p_\alpha : \alpha < \kappa \rangle$, 
where $\kappa = |\mathbb{R}| = |\mathbb{R}^3|$, and constructs by transfinite recursion a sequence of circles $C_\alpha$, each either empty or of radius 1, such that each $C_\alpha$ is disjoint from all previous circles, and if $p_\alpha$ is not contained in a previous 
circle, it is contained in $C_\alpha$. 
Let “large” mean size $\kappa$, and “small” mean size strictly less than $\kappa$. 
The proof turns on the fact that if $\alpha$ is an ordinal $<\kappa$, then $|\alpha|$ is small, and if $\mu, \lambda$ are small cardinals, then $\mu \cdot \lambda$ is small (note this includes the case where both $\mu, \lambda$ are finite).
\footnote{More precisely,  a circle is a locus of points, all lying in a plane and equidistant from some point called the center. 
The proof builds by transfinite recursion a sequence of circles $C_\alpha$, each either empty or of radius 1, such that each 
$C_\alpha$ is disjoint from all previous circles and such that $p_\alpha$ is contained in one of the circles up to and including $C_\alpha$. Suppose we arrive to stage $\alpha$. If $p_\alpha$ is already covered, let $C_\alpha$ be empty and go to the next stage. Otherwise, we first choose a plane $P_\alpha$ containing $p_\alpha$ such that no previous nontrivial circle $C_\beta$ lies in $P_\alpha$. We can do this because there are a large number of planes through $p_\alpha$, and the number of circles chosen so far, $|\alpha|$, is small. Of course, some of the previous circles may intersect our plane. So let us call a point of $P_\alpha$ “bad” if it intersects any earlier $C_\beta$. Each $C_\beta$ can be responsible for at most two bad points, so the number of bad points in $P_\alpha$ is small. It would suffice to choose a circle of radius 1 lying in $P_\alpha$, going through $p_\alpha$ and avoiding all bad points. Consider a fixed bad point b. The number of circles of radius 1 in $P_\alpha$ which can go through both $p_\alpha$ and b is at most two. So each bad point b rules out at most two circles. There are a large number of circles of radius 1 through 
$p_\alpha$ in our plane $P_\alpha$, and a small number which are ruled out, and so we must succeed.}

Induction and recursion on the natural numbers rely on two things: that the natural numbers are well ordered, and that as we walk 
increasingly through them, there are always few things behind us compared to what is ahead of us. Cardinals make these properties available to us 
in walks through a much larger landscape.

\br
\section*{IIIIIIIIIIIIII} 
\br

As model theorists, we see a picture perhaps currently less visible in other areas of mathematics:  mathematical information 
observed across the hierarchy of cardinals can be distilled into powerful definitions and theorems which 
themselves, after the fact, may make no reference to these origins. (As written elsewhere \cite{mm-icm}, perhaps model theory
may have a special ability to make translations of this kind.) 

We now sketch such an example from the {monumental} work of Shelah \cite{Sh:a}. 
We will give just a few details, chosen to illustrate how, along the way, certain properties of cardinals are being used. 

A (first-order) \emph{theory} is a set of axioms, which we can identify with the class of structures satisfying it, 
for example: algebraically closed fields of characteristic zero, real closed fields, or dense linear orders without endpoints. 

We begin with \emph{counting types}. Roughly speaking, a type in model theory tells us what can be said about 
one element of a structure in terms of another set of elements.  Given a particular infinite $\lambda$, 
we say a theory is $\lambda$-stable if over any set of parameters of size $\lambda$ in any model of the theory, 
there are no more than $\lambda$ types; otherwise, it is $\lambda$-unstable. 

If  $K$ is any algebraically closed, hence infinite, field of characterisic $p$ and size 
$\lambda$, then an element can interact with $K$ either by being a root of some minimal polynomial with 
coefficients in $K$ (there are $\lambda = |K|$ such polynomials), or by being transcendental over $K$.  So 
the number of types over $K$ is precisely $\lambda$, and we say the theory of algebraically closed fields of fixed characteristic is 
\emph{$\lambda$-stable for every infinite $\lambda$}.  On the other hand, for every infinite 
$\lambda$ there exists a dense linear order without endpoints of size $\lambda$, call it $M$, 
over which there are strictly more than $\lambda$ cuts; thus, the number of types over $M$ is $>\lambda$. 
The theory of dense linear orders without endpoints is \emph{$\lambda$-unstable for every infinite $\lambda$}. 
A priori, there is perhaps little reason to expect such a pure dichotomy across  
all theories. 

One of Shelah's discoveries was that if a theory $T$ is unstable in some 
$\lambda$ with $\lambda = \lambda^{|T|}$, then there must be a single formula 
causing the number of types to be large. The proof goes by writing the type as a product of its restrictions to each formula, of which there are $|T|$-many, and noticing one of the factors must have size $ > \lambda$, as otherwise the $|T|$-fold product of $\lambda$, i.e. $\lambda^{|T|}$, would have size $\lambda$. 
Thus, instability is a local phenomenon, caused by unstable formulas, which 
we discover by using the arithmetic of a particular cardinal. 
 
Shelah proved a dichotomy theorem:

\begin{quotation}
\noindent\emph{Theorem:} 
Any theory $T$ is either
\begin{itemize}
\item[(a)] stable, i.e. $\lambda$-stable for all infinite lambda such that 
$\lambda = \lambda^{|T|}$. 

\item[(b)] unstable, i.e. $\lambda$-unstable for all infinite lambda.
\end{itemize}
\end{quotation}
This dichotomy, which had a foundational impact on modern model theory, 
is powerful in part because Shelah was able to show that there are 
many other characterizations of this partition.
For example, it may be re-framed in terms of 
ranks, dimensions, definability of types, bounds on the orbits of certain automorphism groups, and also completely ``locally'' in terms of unstable formulas. Let us state the last of these. Say a formula $\vp(\bar{x},\bar{y})$ has the order propery if in models of $T$, for arbitrarily 
large finite $n$, there exist $\langle \bar{a}_i : i < n \rangle$ and $\langle \bar{b}_j : j < n \rangle$ such that $\vp[\bar{a}_i, \bar{b}_j]$ holds if and only if $i < j$. 

\begin{quotation}
\noindent\emph{Theorem:} 
A theory $T$ is unstable if and only if it contains a formula with the \emph{order property}.
 \end{quotation}
By the time we arrive to this theorem, the definition of the order property is simply a kind of long asymmetry, and makes no mention of 
infinite cardinals.     Algebraically closed fields do not have asymmetries of this kind; neither do free groups on fixed finite numbers of generators. 
Infinite random graphs and real closed fields do.

\br

\section*{IIIIIIIIIIIIIII} 

\br
What can the infinite tell us about the finite? 

One can answer this question in many ways; we will give a short example which is recent and close to our hearts. 

Stability has been an extremely productive and powerful influence within model theory, explaining surprising 
similarities and
differences among 
classes of infinite structures (such as algebraically closed versus real closed fields). 
It long seemed to be about infinite structures in 
a deep way.  Recently, the picture has shifted. 

Consider an a priori different setting. 
The finite Ramsey theorem tells us that for any finite $n$ there is 
$M=M(n)$ such that any finite graph $G$ of size at least $M$ contains a homogeneous set of size $n$.  
Approximately, a graph of size $n$ will have a homogeneous set of size at least $\log n$.  Preciesely determining the bounds or ``Ramsey numbers'' for all graphs or for certain classes of graphs leads to much challenging combinatorics. 

Homogeneous sets may be thought of as special cases of the following model theoretic definition, which uses all formulas, not only the formula expressing whether an edge holds. (This connection was already made in Ramsey’s original paper \cite{ramsey}.) 

\emph{Definition.} In any mathematical structure, call a sequence of elements indiscernible if for any finite k, any two strictly increasing subsequences of length k satisfy the same first order formulas.

Part of the work on stability in \cite{Sh:a} 
is the 
following surprising theorem about infinite models of stable theories: 
\emph{If T is stable then in any infinite model M of T of size $\kappa^+$ there is a set of indiscernibles of size $\kappa^+$.}
Informally, in infinite models of stable theories of reasonable sizes, one can extract indiscernible sequences (indeed, sets) of the same size. 
As an example, in any algebraically closed field of fixed characteristic and size $\aleph_1$, one can find a transcendence basis of size $\aleph_1$ (whose elements in some sense interact completely generically with each other).

It turns out that this theorem about infinite stable theories makes a precise and useful prediction about the finite world. 
A 2014 theorem of Malliaris and Shelah [9] shows that one can in some sense “finitize” stability in the case of graphs by forbidding a certain family of configurations, the half-graphs, which are finite versions of the order property for the graph edge relation. It becomes possible to build an approximation to the structural understanding of stability within these ``stable-like'' graphs, and thus to prove, mirroring the infinite case, that in so-called finite stable graphs there are much larger homogeneous sets than predicted by Ramsey theory. A simple version of the 
resulting theorem says: 
\emph{There is a constant $c = c(k)$ such that if $G$ is a finite graph which does not contain any $k$-half graph then $G$
homogeneous set of size $|G|^c$.} 

In this simple example, and in analogous refinements of Szemer\'edi regularity \cite{MiSh:978} \cite{MiSh:E84}, one uses an understanding of infinite stability phenomena (themselves revealed in the infinite world through an understanding of the hierarchy of cardinals) as a precise guide to prove theorems about the finite world.

\br
\section*{IIIIIIIIIIIIIIII} 
\br

If in the course of this paper we have invoked 
Definition \ref{d:defn} to start to see how, in their first century and half, 
Cantor's infinite cardinals show themselves to be delicately and deeply interwoven with the fabric of our mathematical world and our mathematical thought,  
Definition \ref{d:defn} shouldn't however be used \emph{against} cardinals, as a way of silencing them 
or banishing them from hypotheses of theorems. 

There is a greatness in certain mathematical ideas -- as also in certain works of fiction -- which 
challenges us 
by enlarging the boundaries of what we thought possible without ever directly breaking with the 
world we knew.  Emerson, writing on character, says of 
certain great figures that ``the largest part of their 
power was latent.''  Mathematics, having as it does a longer history, 
may take centuries or even millenia to understand and to assimilate a new definition, but 
this doesn't mean its power isn't felt already in our time.


\vspace{15mm}

\begin{thebibliography}{50}

\bibitem{cohen} P. Cohen, \emph{Set Theory and the Continuum Hypothesis}, Dover edition, 2008. Page 2.  

\bibitem{emerson} R. W. Emerson.  \emph{Essential writings of Ralph Waldo Emerson.} Modern Library. Page 327.

\bibitem{gowers}  
W. T. Gowers, ``Lower bounds of tower type for Szemer\'edi's uniformity lemma.'' Geometric and Functional Analysis, vol. 7 (1997) 322--337.

\bibitem{rs} A. Grothendieck, \emph{Recoltes et Semailles}. 

\bibitem{kt} Komjath and Totik. \emph{Problems and Theorems in Classical Set Theory.}  Problem Books in Mathematics, Springer, 2006.  
Section 13. 

\bibitem{k-s} J. Koml\'os and M. Simonovits, ``Szemer\'edi's Regularity Lemma and its applications in graph theory.'' (1996) In \emph{Combinatorics: Paul Erd\"os is Eighty}, Vol. 2 (D. Mikl\'os, V. T. S\'os and T. Sz\"onyi, eds), Bolyai Society Math. Studies, Keszthely, Hungary, pp. 295-352.

\bibitem{mm-icm} M. Malliaris. ``Model theory and ultraproducts.''  
\emph{Proceedings of the International Congress of Mathematicians—Rio de Janeiro 2018}. 
Vol. II. Invited lectures, 83--97, World Sci. Publ., Hackensack, NJ, 2018.

\bibitem{mm-ap1} M. Malliaris and A. Peretz, ``What simplicity is not.'' 
\emph{Simplicity: ideals of practice in mathematics and the arts}. 51--58, Math. Cult. Arts, Springer, Cham, 2017.

\bibitem{MiSh:978} M. Malliaris and S. Shelah. ``Regularity lemmas for stable graphs.''  Trans. Amer. Math Soc, 366 (2014), 1551--1585.

\bibitem{MiSh:E84} M. Malliaris and S. Shelah. ``Notes on stable regularity.''   Bull. Symb. Log. 27 (2021), no. 4, 415--425. 

\bibitem{ap} A. Peretz,  ``Geometry of forking in simple theories.'' J Symbolic Logic 71, 1 (2006), 347--359.  

\bibitem{ramsey} F. P. Ramsey, ``On a problem of formal logic.'' Proc London Math Soc 30 (1929) 338--384.

\bibitem{Sh:a} 
S. Shelah, \emph{Classification Theory and the number of non-isomorphic models}, first edition, 1978. (rev. ed. 1990)  North-Holland. 

\end{thebibliography}
\end{document}